\documentclass{amsart}

\author{Samuel Boissi{\`e}re \and Alessandra Sarti}

\title{Contraction of excess fibres between the McKay correspondences in dimensions two and three}

\date{Monday, april 15, 2005}

\address{Samuel Boissi{\`e}re, Fachbereich f{\"u}r Mathematik, Johannes Gutenberg-Universit{\"a}t,
55099 Mainz, Germany}

\email{boissiere@mathematik.uni-mainz.de}

\urladdr{http://sokrates.mathematik.uni-mainz.de/$\sim$samuel}

\address{Alessandra Sarti, Fachbereich f{\"u}r Mathematik, Johannes Gutenberg-Universit{\"a}t,
55099 Mainz, Germany}

\email{sarti@mathematik.uni-mainz.de}

\urladdr{http://www.mathematik.uni-mainz.de/$\sim$sarti}

\usepackage[english]{babel}
\usepackage{amsmath,amsfonts,amssymb}
\usepackage{bbm,dsfont}
\usepackage{mathrsfs}
\usepackage[all]{xy}
\usepackage{ifthen}
\usepackage{rotating}

\DeclareMathOperator{\Irr}{Irr} \DeclareMathOperator{\Id}{Id}
\DeclareMathOperator{\rank}{rk} \DeclareMathOperator{\GL}{GL}
\DeclareMathOperator{\SO}{SO} \DeclareMathOperator{\SU}{SU}
\DeclareMathOperator{\Sym}{S} \DeclareMathOperator{\SL}{SL}
\DeclareMathOperator{\Jac}{Jac} \DeclareMathOperator{\Spec}{Spec}

\newcommand{\ie}{\textit{i.e. }}
\newcommand{\loccit}{\textit{loc.cit. }}

\newcommand{\Hilbf}[3]{{#1}\ifthenelse{\equal{#1}{}}{}{\text{-}}\mathcal{H}ilb^{#2}_{#3}}
\newcommand{\Hilb}[3]{{#1}\ifthenelse{\equal{#1}{}}{}{\text{-}}\mathrm{Hilb}^{#2}{\left(#3\right)}}

\newcommand{\ci}{\mathrm{i}}

\newcommand{\IC}{\mathds{C}}
\newcommand{\IH}{\mathds{H}}
\newcommand{\IR}{\mathds{R}}
\newcommand{\IZ}{\mathds{Z}}

\newcommand{\IP}{\mathbb{P}}
\newcommand{\IS}{\mathbb{S}}

\newcommand{\Iun}{\mathbbm{1}}
\newcommand{\Ii}{\mathbbm{i}}
\newcommand{\Ij}{\mathbbm{j}}
\newcommand{\Ik}{\mathbbm{k}}

\newcommand{\cI}{\mathcal{I}}
\newcommand{\cO}{\mathcal{O}}
\newcommand{\cQ}{\mathcal{Q}}
\newcommand{\cT}{\mathcal{T}}
\newcommand{\cX}{\mathcal{X}}
\newcommand{\cY}{\mathcal{Y}}
\newcommand{\cZ}{\mathcal{Z}}

\newcommand{\cH}{\mathscr{H}}
\newcommand{\cS}{\mathscr{S}}

\newcommand{\fA}{\mathfrak{A}}
\newcommand{\fG}{\mathfrak{G}}
\newcommand{\fI}{\mathfrak{I}}
\newcommand{\fR}{\mathfrak{R}}
\newcommand{\fS}{\mathfrak{S}}

\newcommand{\fm}{\mathfrak{m}}
\newcommand{\fn}{\mathfrak{n}}

\newtheorem{theorem}{Theorem}[section]
\newtheorem{lemma}[theorem]{Lemma}
\newtheorem{proposition}[theorem]{Proposition}
\newtheorem{corollary}[theorem]{Corollary}
\newtheorem{remark}[theorem]{Remark}

\begin{document}

\begin{abstract}
The quotient singularities of dimensions two and three obtained from polyhedral
groups and the corresponding binary polyhedral groups admit natural resolutions
of singularities as Hilbert schemes of regular orbits whose exceptional fibres
over the origin reveal similar properties. We construct a morphism between
these two resolutions, contracting exactly the excess part of the exceptional
fibre. This construction is motivated by the study of some pencils of K3-surfaces
arising as minimal resolutions of quotients of nodal surfaces with high symmetries. 
\end{abstract}



\maketitle

\pagestyle{myheadings}

\markboth{SAMUEL BOISSI{\`E}RE AND ALESSANDRA SARTI}{Contraction of excess
fibres}

\section{Introduction}

Consider a binary polyhedral group $\widetilde{G}\subset \SU(2)$
corresponding to a polyhedral group $G\subset \SO(3,\IR)$ through the
double-covering $\SU(2)\rightarrow \SO(3,\IR)$. The group $\widetilde{G}$
acts freely on $\IC^2$ and the quotient $\IC^2/\widetilde{G}$ is a surface
singularity with an isolated singular point at the origin. The exceptional divisor
of its minimal resolution of singularities $\cX\rightarrow \IC^2/\widetilde{G}$
is a bunch of smooth rational curves of self-intersection $-2$, intersecting
transversely, whose intersection graph is an A-D-E Dynkin diagram. The classical
McKay correspondence (\cite{McKay}) relates this intersection graph to the
representations of the group $\widetilde{G}$, associating bijectively each
exceptional curve to a non-trivial irreducible representation of the group: the
correspondence in fact identifies the intersection graph
with the McKay quiver of the action of $\widetilde{G}$ on $\IC^2$. Among these
irreducible representations we find all irreducible representations of the
group $G$: we call them \emph{pure} and the remaining ones \emph{binary}. Since
$\widetilde{G}/G\cong \{\pm 1\}$, one can produce a $G$-invariant cone
$\IC^2/\{\pm 1\}\xrightarrow{\sim}K\hookrightarrow \IC^3$ whose quotient $K/G$
is isomorphic to $\IC^2/\widetilde{G}$. In this note, we prove the following
result, conjectured by W.~ P.~ Barth:

\begin{theorem}
\label{th:barth} There exists a crepant resolution of singularities of
$\IC^3/G$ containing a partial resolution $\cY\rightarrow K/G$ with the
property that the intersection graph of its exceptional locus is precisely the
McKay quiver of the action of $G$ on $\IC^3$, together with a resolution map
$\cX\rightarrow \cY$ mapping isomorphically the exceptional curves
corresponding to pure representations and contracting those associated with
binary representations to ordinary nodes.
\end{theorem}

We make this construction in the framework of the Hilbert schemes of regular
orbits of Nakamura (\cite{Nakm}) providing, thanks to the Bridgeland-King-Reid
theorem (\cite{BKR}), the natural candidates for the resolutions of
singularities in dimensions two and three. We produce a morphism $\cS$ between
these two resolutions of singularities, define our partial resolution $\cY$ as
the image of this map and study the effect of $\cS$ on the exceptional
fibres:
$$
\xymatrix{\Hilb{\widetilde{G}}{}{\IC^2}\ar[rr]^\cS\ar[dd]_{\widetilde{\pi}}\ar@{->>}[dr]_\cS
&&\Hilb{G}{}{\IC^3}\ar[dd]^\pi \\
& \cY\ar[d]\ar@{^(->}[ru] & \\
\IC^2/\widetilde{G}\ar[r]^{\sim} &K/G\ar@{^(->}[r]& \IC^3/G}
$$
Although the exceptional fibres can be described very explicitly in all cases
(see \cite{INm1}), by principle our proof avoids any case-by-case analysis.
Therefore, the key point consists in a systematic modular interpretation of the objects at issue.

From the point of view of the story of the McKay correspondence, this construction
shows some news properties revealing again the fertility of the geometric
construction of the McKay correspondence following Gonzales-Sprinberg and
Verdier \cite{GSV}, Ito-Nakamura \cite{INm1}, Ito-Najakima \cite{INj} and Reid
\cite{Reid1}. The beginning of the story was devoted to the study of all
situations in dimension two and three, in general by a case-by-case analysis.
Then efforts were made to understand how to get all these cases by one general
geometric construction (\cite{INj,BKR}). The development followed then the
cohomological direction in great dimensions in a symplectic setup (\cite{K1,GK}), leading to
an explicit study of a family of examples of increasing dimension
for the specific symmetric group problem (\cite{Boi}). The new point of view in the present
paper consists in working between two situations of different dimensions for
different - but related - groups and construct a relation between them. This may
be considered as a concrete application of some significant results in this
area coming again at the beginning of the story, dealing with a now quite classical
material approached by natural transformations between moduli spaces.

This study is motivated by previous works of Sarti \cite{Sa} and Barth-Sarti \cite{BaSa} studying
special pencils of surfaces in $\IP_3$ with bipolyhedral symmetries. The minimal resolutions of the 
associated quotient surfaces are K3-surfaces with maximal Picard numbers. For some special fibres
of these pencils, the resolution looks locally like the quotient of a cone by a polyhedral group, and
our result gives a local interpretation of the exceptional locus in these cases.

The structure of the paper is as follows: in Section \ref{s:clusters} we introduce the notations and
 we recall some basic facts about clusters and in Section \ref{moduli} we recall the construction of the Hilbert schemes
 of points and clusters. The Sections \ref{binary}, \ref{graph} and \ref{exceptional} give a brief survey
 on polyhedral, binary polyhedral and bipolyhedral groups, their representations and the classical Mckay correspondences in dimensions two and three. 
In Section \ref{s:geomconstr} we start the study of the map $\cS$. First we show that it is well defined 
(lemma \ref{bendef}) and then that it is a regular projective map, which induces a map between the exceptional 
fibres (proposition \ref{regpro}). In Section \ref{sezteoremone} the theorem \ref{th:theoremone} is the 
fundamental step for proving the main theorem \ref{th:barth}: we show that the map $\cS$ contracts the
 curves corresponding to the binary representations and maps the curves corresponding to the pure 
representations isomorphically to the exceptional curves downstairs. In Section \ref{s:examplecyclic}, 
as an example we describe in details the case when $\widetilde{G}$ is a cyclic group. Finally the Section 
\ref{ultimo} is devoted to an application to resolutions of pencils of K3-surfaces. 

\smallskip

{\bf Acknowledgments:} we thank Wolf Barth for suggesting the problem and for many helpful comments, Manfred Lehn for his invaluable help during the preparation of this paper and Hiraku Nakajima for interesting explanations.

\section{Clusters}
\label{s:clusters}

In the sequel, we aim to study a link between the two-dimensional and the
three-dimensional McKay correspondences. In order to avoid confusion, we shall
use different sets of letters for the corresponding algebraic objects at issue
in both situations. In this section, we fix the notations and the terminology.

\subsection{General setup}
\label{ss:generalsetup}

Let $V$ be a $n$-dimensional complex vector space and $\fG$ a finite subgroup
of $\SL(V)$. We denote by $\cO(V):=\Sym^*(V^\vee)$ the algebra of  polynomial
functions on $V$, with the induced left action $g\cdot f:=f\circ g^{-1}$ for
$f\in \cO(V)$ and  $g\in \fG$.

We choose a basis $X_1,\ldots,X_n$ of linear forms on $V$, denote the ring of
polynomials in $n$ indeterminates by $S:=\IC[X_1,\ldots,X_n]$ and identify
$\cO(V)\cong S$. The ring $S$ is given a  graduation by the total degree of a
polynomial, where each indeterminate $X_i$ has degree $1$. In particular, the
action of the group $\fG$ on $S$ preserves the degree.

Let $\fm_S:=\langle X_1,\ldots,X_n\rangle$ be the maximal ideal of $S$ at the
origin. We denote by $S^\fG$ the subring of $\fG$-invariant polynomials, by
$\fm_{S^\fG}$ its maximal ideal at the origin and by $\fn_\fG:=\fm_{S^\fG}\cdot
S$ the ideal of $S$ generated by the non-constant $\fG$-invariant polynomials
vanishing at the origin. The quotient ring of \emph{coinvariants} is by
definition $S_\fG:=\left. S\right/\fn_\fG$.

An ideal $\fI\subset S$ is called a $\fG$-\emph{cluster} if it is globally
invariant under the action of $\fG$ and the quotient
$S/\fI$ is isomorphic, as a $\fG$-module, to the regular representation of
$\fG$: $S/\fI\cong \IC[\fG]$. A closed subscheme $Z\subset \IC^n$ is called a
$\fG$-\emph{cluster} if its defining ideal $\fI(Z)$ is a $\fG$-cluster. Such a
subscheme is then zero-dimensional and has length $|\fG|$. For instance, a free
$\fG$-orbit defines a $\fG$-cluster. In particular, a $\fG$-cluster contains
only one orbit: the support of a cluster is a union of orbits, and any
function constant on one orbit and vanishing on another one would induce a different copy
of the trivial representation in the quotient $S/\fI$.

We are particularly interested in $\fG$-clusters supported at the origin. Then
$\fI\subset~ \fm_S$ and in fact this condition is enough to assert that the
cluster is supported at the origin: else, the support of the cluster would
consist in more than one orbit.
Furthermore, one has automatically $\fn_\fG\subset \fI$, since any non-constant
function $f\in \fn_\fG$ not contained in $\fI$ would induce a new copy of the
trivial representation in the quotient $S/\fI$, different from the one already given by
the constant functions. Hence we wish to understand the structure of the $\fG$-clusters
$\fI$ such that $\fn_\fG\subset \fI\subset \fm_S$, equivalent to the
study of the quotient ideals $\fI/\fn_\fG \subset \fm_S/\fn_\fG \subset
S/\fn_\fG=S_\fG$, with the exact sequence:
\begin{equation}\label{exactseq:cluster}
0\longrightarrow \fI/\fn_\fG \longrightarrow S_\fG \longrightarrow S/\fI \longrightarrow 0.
\end{equation}

From now on, we assume that the group $\fG$ is a subgroup of index $2$ of a
group  $\fR\in~ \GL(V)$ generated by \emph{reflections} (we follow here the terminology of \cite{Co}),
\ie elements $g\in \fR$ such that $\rank (g-~\Id_V)=~1$.

The structure of the action of $\fR$ on $S$ has the following properties (see \cite{Co}):
\begin{itemize}
\item The algebra of invariants $S^{\fR}$ is a polynomial algebra generated by
exactly $n$ algebraically independent homogeneous polynomials $f_1,\ldots,f_n$
of degrees $d_i$. \item $|\fR|=d_1\cdot\ldots\cdot d_n$. \item The set of
degrees $\{d_1,\ldots,d_n\}$ is independent of the choice of the homogeneous
generators. \item The algebra of coinvariants is isomorphic to the regular
representation: $S_{\fR}\cong \IC[\fR]$.
\end{itemize}
As a byproduct, we get that the algebra of coinvariants $S_{\fR}$ is a graded
finite-dimensional algebra.

From this and the fact that $\fG=\fR\cap \SL(V)$, one deduces the structure  of
the action of $\fG$ on $S$ (see \cite{Bou,GNS1,GNS2}):
\begin{itemize}
\item There exists a homogeneous $\fR$-skew-invariant polynomial $f_{n+1}\in
S$, \ie  such that $g\cdot f_{n+1}=\det(g).f_{n+1}$ for all $g\in \fR$, unique
up to a multiplicative constant, dividing any $\fR$-skew-invariant polynomial:
hence the set $f_{n+1}\cdot S^{\fR}$ is precisely the set of
$\fR$-skew-invariants. A natural choice for this element is
$f_{n+1}=\Jac(f_1,\ldots,f_{n})$. \item $S^\fG=\IC[f_1,\ldots,f_n,f_{n+1}]$.
\item $\fn_{\fG}=\fn_{\fR}\oplus \IC f_{n+1}$. \item $S_{\fR}=S_{\fG}\oplus
\IC f_{n+1}$.
\end{itemize}
Note that, as a $\fG$-module, $\IC[\fR]$ is isomorphic to two copies of
$\IC[\fG]$.  It follows that $\fm_S/\fn_\fG$ is a graded finite-dimensional
algebra which, as a $\fG$-module, consists exactly of each non-trivial
representation $\rho$ of $\fG$ repeated $2\dim \rho$ times: one can denote the
occurrences of each representation $\rho$ by $V^{(1)}(\rho),\ldots,V^{(2\dim
\rho)}(\rho)$ where each $V^{(i)}(\rho)$ is given by homogeneous polynomials
modulo $\fn_\fG$.

Thanks to the exact sequence (\ref{exactseq:cluster}), giving a $\fG$-cluster
supported at the origin consists in choosing, for each non-trivial
representation $\rho$ of $\fG$, $\dim \rho$ copies of $\rho$ in
$\fm_S/\fn_\fG$. But this gives many choices since any linear combination of
some $V^{(i)}(\rho)$ and $V^{(j)}(\rho)$ provides such a copy. The ground idea is that one does not
have to make all these choices in order to define $\fI$ (see \S\ref{s:examplecyclic} for an explicit example).

For such an ideal $\fI$ with $\fn_\fG\subset \fI\subset \fm_S$, we consider the
finite-dimensional $\fG$-modules $W\subset S$ generating $\fI$ in the sense
that $\fI=W\cdot S+\fn_\fG$. Such modules do exist thanks to the preceding
construction. Among these choices, we consider the minimal ones, \ie such that
no strict $\fG$-submodule of them generate $\fI$ in the preceding sense.

If $W$ is a generator in this sense, then 
$$
\fI=W\cdot S+\fn_\fG=W+\fm_S\cdot W+\fn_\fG=W+\fm_S\cdot \fI+\fn_\fG.
$$
This means that the $\IC$-linear map $W\rightarrow \fI/(\fm_S\cdot
\fI+\fn_\fG)$  is surjective. Also, since $W$ is a $\fG$-module and since
$\fm_S\cdot \fI+\fn_\fG$ is $\fG$-stable, this map is $\fG$-linear. If $W$ is a
minimal set of generators,  it satisfies in particular $W\cap (\fm_S\cdot
\fI+\fn_\fG)=\{0\}$ since this intersection would provide a $\fG$-submodule
whose complementary in $W$ is a smaller $\fG$-submodule generating $\fI$.
Hence, for $W$ minimal one gets an isomorphism of $\fG$-modules $W\cong
\fI/(\fm_S\cdot \fI+\fn_\fG)$. We set then $V(\fI):=\fI/(\fm_S\cdot
\fI+\fn_\fG)$. The set of generators of $V(\fI)$ may not be uniquely
determined, but its structure as a $\fG$-module is unique. The important issue,
that will be the core of the classification, will be to determine whether
$V(\fI)$ is irreducible or not, although it is a minimal set of generators.

\subsection{Notations for the two- and three-dimensional cases}

When applying the preceding constructions in dimensions two or three, we fix
the  following notations:
\begin{itemize}
\item For $n=2$, the polynomial ring is denoted by $A:=\IC[x,y]$, the group by
$\widetilde{G}$ and any ideal by $I$.

\item For $n=3$, the polynomial ring is denoted by $B:=\IC[a,b,c]$, the group
by $G$ and any ideal by $J$.
\end{itemize}

\section{Moduli space of clusters}\label{moduli}

We recall here the constructions of the Hilbert schemes of points or clusters.

\subsection{Hilbert scheme of points}

Let $X\subset \IP^n_\IC$ be a quasi-projective scheme and $N$ a positive integer. Consider the contravariant
functor $\Hilbf{}{N}{X}$ from the category of schemes to the category of sets
$$
\Hilbf{}{N}{X}:(Schemes)\rightarrow (Sets)
$$
which is given by
$$
\Hilbf{}{N}{X}(T):=
\left\{ Z\subset T\times X\left|
\begin{array}{ll}
(a) & Z \text{ is a closed subscheme}\\
(b) & \text{the morphism } Z\hookrightarrow T\times X\xrightarrow{p}T\text{ is flat}\\
(c) & \forall t\in T, Z_t\subset X \text{ is a closed subscheme}\\
& \text{of dimension $0$ and length $N$}
\end{array}
\right. \right\}
$$

By a theorem of Grothendieck (\cite{Gro}), this functor is representable by a
quasi-projective scheme $\Hilb{}{N}{X}$ equipped with a \emph{universal family} $\Xi_N^X\subset
\Hilb{}{N}{X}\times X$. In the sequel, we shall always denote by $p$ the
projection to the moduli space (here $\Hilb{}{N}{X}$) and by $q$ the projection
to the base (here $X$). When $X$ is projective, the scheme $\Hilb{}{N}{X}$ is
projective and comes with a very ample line bundle (for $\ell\gg 0$):
$$
\det\left(p_*\left(\cO_{\Xi_N^X}\otimes q^*\cO_X(\ell)\right)\right).
$$
When $X=\IC^n$, one gets an open immersion
$\Hilb{}{N}{\IC^n}\hookrightarrow\Hilb{}{N}{\IP^n_\IC}$ corresponding to the
restriction of the universal family. The induced restriction of the preceding
determinant line bundle provides us the very ample line bundle 
$\det\left(p_*\cO_{\Xi_N^{\IC^n}}\right)$ on $\Hilb{}{N}{\IC^n}$.

There exists a natural projective morphism from $\Hilb{}{N}{X}$ to the
symmetric  product $\Sym^{N}(X)$ sending a closed subscheme to the
corresponding $0$-cycle describing its support, called the \emph{Hilbert-Chow}
morphism:
$$
\cH:\Hilb{}{N}{X}\longrightarrow \Sym^{N}(X).
$$
By a theorem of Fogarty (\cite{Fo}), the scheme $\Hilb{}{N}{X}$ is connected.
For $\dim X=2$, it is reduced, smooth and the morphism $\cH$ is a resolution of
singularities.

\subsection{Hilbert scheme of regular orbits}

We consider the sub-functor $\Hilbf{\fG}{}{\IC^n}$ of $\Hilbf{}{|\fG|}{\IC^n}$
given by
$$
\Hilbf{\fG}{}{\IC^n}(T):=\left\{ Z\in \Hilbf{}{|\fG|}{\IC^n}(T)\,|\, \forall t\in T, Z_t\subset \IC^n \text{ is a $\fG$-cluster}\right\}.
$$
This functor is representable by a quasi-projective scheme
$\Hilb{\fG}{}{\IC^n}$ called the
\emph{Hilbert scheme of $\fG$-regular orbits}, which is a union of some
connected components of the subscheme of $\fG$-fixed points
$\left(\Hilb{}{|\fG|}{\IC^n}\right)^\fG$. Furthermore, the quotient $\IC^n/\fG$
can be identified with a closed subscheme of $\Sym^{|\fG|}(\IC^n)$ and since the
support of a $\fG$-cluster consists exactly of one orbit through $\fG$, the
restriction of the Hilbert-Chow morphism factorizes through a projective
morphism (see \cite{BKR,INj,T}):
$$
\cH:\Hilb{\fG}{}{\IC^n}\longrightarrow \IC^n/\fG.
$$
There is a unique irreducible component of
$\Hilb{\fG}{}{\IC^n}$ containing the free $\fG$-orbits and mapping birationally
onto $\IC^n/\fG$. This component is taken as the definition of the Hilbert
scheme of $\fG$-regular orbits in \cite{Nakm}. By the theorem of
Bridgeland-King-Reid \cite{BKR}, if $n\leq 3$, then $\Hilb{\fG}{}{\IC^n}$ is
already irreducible, reduced, smooth and the map $\cH$ a crepant resolution of
singularities of the quotient $\IC^n/\fG$. Moreover, $\cH$ is an isomorphism
over the open subset of free $\fG$-orbits. As a byproduct, the two definitions
coincide.

As before, the scheme $\Hilb{\fG}{}{\IC^n}$ is equipped with a universal family
$\cZ_\fG$ which is the restriction of the universal family
$\Xi_{|\fG|}^{\IC^n}$ corresponding to the closed immersion
$\Hilb{\fG}{}{\IC^n}\hookrightarrow \Hilb{}{|\fG|}{\IC^n}$. The induced
restriction of the determinant line bundle provides us, by naturality of the
construction of the determinant of a family (see \cite[\S 8.1]{HL}), the very
ample line bundle $\det\left(p_*\cO_{\cZ_\fG}\right)$ on $\Hilb{\fG}{}{\IC^n}$.

\section{Rotation groups}
\label{binary}

\subsection{Polyhedral groups}

Let $\SO(3,\IR)$ be the group of rotations in $\IR^3$. Up to conjugation,
there are five different types of finite subgroups of $\SO(3,\IR)$, called
\emph{polyhedral groups}:
\begin{itemize}
\item the cyclic groups $C_n\cong \IZ/n\IZ$ of order $n\geq 1$;

\item the dihedral groups $D_n\cong \IZ/n\IZ \rtimes \IZ/2\IZ$ of order $2n$,
$n\geq 1$;

\item the group $\cT$ of positive isometries of a regular tetrahedra,
isomorphic to the alternate group $\fA_4$ of order $12$;

\item the group $\cO$ of positive isometries of a regular octahedra or a cube,
isomorphic to the symmetric group $\fS_4$ of order $24$;

\item the group $\cI$ of positive isometries of a regular icosahedra or a
regular dodecahedra, isomorphic to the alternate group $\fA_5$ of order $60$.
\end{itemize}

\subsection{Binary polyhedral groups}

Let $\IH$ be the real algebra of quaternions, with basis $(\Iun,\Ii,\Ij,\Ik)$. 
The \emph{norm} of a quaternion $q=a\cdot\Iun+b\cdot\Ii+c\cdot\Ij+d\cdot\Ik$ is 
$N(q):=~a^2+~b^2+~c^2+~d^2$, $a,b,c,d\in \IR$. Let $\IS$ be the three-dimensional sphere of
quaternions of length $1$ and $H$ the three-dimensional vector subspace of \emph{pure} quaternions 
(\ie $a=0$). For $q\in \IS$, the action by conjugation $\phi(q):~H\rightarrow~ H$, $x\mapsto q\cdot x\cdot q^{-1}$ is an isometry. Since the group $\IS$ is isomorphic to $\SU(2)$ by the identification
$$
q=\left(
\begin{matrix}
a+\ci b& c+\ci d\\
-c+\ci d& a-\ci b
\end{matrix}
\right),
$$ 
one gets an exact sequence
$$
0\longrightarrow \{\pm 1\} \longrightarrow \SU(2)
\overset{\phi}{\longrightarrow} \SO(3,\IR)\longrightarrow 0.
$$

For any finite subgroup $G\subset\SO(3,\IR)$, the inverse image
$\widetilde{G}:=\phi^{-1}G$ is called a \emph{binary polyhedral group}. It is a
finite  subgroup of $\SU(2)$ or equivalently, up to conjugation, of
$\SL(2,\IC)$:
\begin{itemize}
\item the binary cyclic groups $\widetilde{C}_n\cong C_{2n}$ have order
$2n$;

\item the binary dihedral groups $\widetilde{D}_n$ have order $4n$;

\item the binary tetrahedral group $\widetilde{\cT}$ has order $24$;

\item the binary octahedral group $\widetilde{\cO}$ has order $48$;

\item  the binary icosahedral group $\widetilde{\cI}$ has order $120$.
\end{itemize}

\subsection{Representations of polyhedral groups}
\label{ss:representations}

Consider a binary polyhedral group $\widetilde{G}$, the associated polyhedral group
$G$ and set $\tau:=\{\pm 1\}$:
$$
0\longrightarrow \tau \longrightarrow \widetilde{G} \overset{\phi}{\longrightarrow}
G\longrightarrow 0.
$$

This exact sequence induces an  injection of the set of irreducible
representations of $G$ in the set of irreducible representations of $\widetilde{G}$: if
$\rho:G\rightarrow \GL(V)$ is an irreducible representation of $G$, it induces
by composition a representation of $\widetilde{G}$ which is $\tau$-invariant,
\ie such that $\rho(-g)=\rho(g)$ for all $g\in \widetilde{G}$. Thanks to this property,
if the representation $\rho$ would admit a non-trivial $\widetilde{G}$-submodule, it
would also be a non-trivial $G$-submodule after going to the quotient
$\widetilde{G}/\tau\cong G$. This shows also that the image of the
injection (since $G$ is a quotient of $\widetilde{G}$):
$$
\Irr(G)\hookrightarrow \Irr(\widetilde{G})
$$
consists precisely on those irreducible representations which are
$\tau$-invariant. These representations are called \emph{pure} and the remaining
representations are called \emph{binary}. More precisely, if
$\rho:\widetilde{G}\rightarrow \GL(V)$ is an irreducible representation of
$\widetilde{G}$, the subspace
$$
V^\tau:=\{v\in V\,|\,v=\rho(-1)v\}
$$
is a $\widetilde{G}$-submodule of $V$. Hence either $V^\tau=V$ and the
representation  $\rho$ is pure, or $\rho$ is binary and $V^\tau=\{0\}$.

For each type of binary polyhedral group, we draw the list of the irreducible
representations with their dimension. The binary representations are labelled
by a ``$\,\widetilde{\text{ }}\,$'' and the trivial representation is denoted
by  $\chi_0$ in all cases:

\begin{itemize}
\item binary cyclic group $\widetilde{C}_n$, $n\geq 1$:

\begin{tabular}{ |c||c|c|c|}
  \hline
   representation & $\chi_0$ & $\{\chi_j\}_{j=1,\ldots,n-1}$ & $\{\widetilde{\chi}_j\}_{j=1,\ldots,n}$ \\
  \hline
  dimension & 1 & 1 & 1 \\
  \hline
\end{tabular}

\item binary dihedral group $\widetilde{D}_n$ for $n=2\ell+1$, $\ell\geq 1$:

\begin{tabular}{ |c||c|c|c|c|c|c|}
  \hline
   representation & $\chi_0$ & $\chi_1$ & $\{\tau_j\}_{j=1,\ldots,\ell}$ &
    $\widetilde{\chi}_1$& $\widetilde{\chi}_2$ & $\{\widetilde{\sigma}_j\}_{j=1,\ldots,\ell}$ \\
  \hline
  dimension & 1 & 1 & 2 & 1 & 1 & 2 \\
  \hline
\end{tabular}

\item binary dihedral group $\widetilde{D}_n$ for $n=2\ell$, $\ell\geq 1$:

\begin{tabular}{ |c||c|c|c|c|c|c|}
  \hline
   representation & $\chi_0$ & $\chi_1$ & $\chi_2$ & $\chi_3$ & $\{\tau_j\}_{j=1,\ldots,\ell-1}$ &
     $\{\widetilde{\sigma}_j\}_{j=1,\ldots,\ell}$ \\
  \hline
  dimension & 1 & 1 & 1 & 1 & 2 & 2 \\
  \hline
\end{tabular}

\item binary tetrahedral group $\widetilde{\cT}$:

\begin{tabular}{ |c||c|c|c|c|c|c|c|}
  \hline
   representation & $\chi_0$ & $\chi_1$ & $\chi_2$ & $\chi_3$ &
   $\widetilde{\chi}_1$ & $\widetilde{\chi}_2$ & $\widetilde{\chi}_3$ \\
  \hline
  dimension & 1 & 1 & 1 & 3 & 2 & 2 & 2\\
  \hline
\end{tabular}

\item binary octahedral group $\widetilde{\cO}$:

\begin{tabular}{ |c||c|c|c|c|c|c|c|c|}
  \hline
   representation & $\chi_0$ & $\chi_1$ & $\chi_2$ & $\chi_3$ & $\chi_4$ &
   $\widetilde{\chi}_1$ & $\widetilde{\chi}_2$ & $\widetilde{\chi}_3$ \\
  \hline
  dimension & 1 & 1 & 2 & 3 & 3 & 2 & 2 & 4\\
  \hline
\end{tabular}

\item binary icosahedral group $\widetilde{\cI}$:

\begin{tabular}{ |c||c|c|c|c|c|c|c|c|c|}
  \hline
   representation & $\chi_0$ & $\chi_1$ & $\chi_2$ & $\chi_3$ & $\chi_4$ &
   $\widetilde{\chi}_1$ & $\widetilde{\chi}_2$ & $\widetilde{\chi}_3$ & $\widetilde{\chi}_4$ \\
  \hline
  dimension & 1 & 3 & 3 & 4 & 5 & 2 & 2 & 4 & 6\\
  \hline
\end{tabular}
\end{itemize}

\subsection{Bipolyhedral groups}

For $p,q\in \IS$, the action $\sigma(p,q):\IH\rightarrow \IH, x\mapsto p\cdot x\cdot q^{-1}$ is an isometry and one
gets an exact sequence
$$
0\longrightarrow \{\pm 1\} \longrightarrow \SU(2)\times\SU(2)
\overset{\sigma}{\longrightarrow} \SO(4,\IR)\longrightarrow 0.
$$
For any binary polyhedral group $\widetilde{G}$,
the direct image $\sigma(\widetilde{G}\times\widetilde{G})\subset \SO(4,\IR)$ is called a \emph{bipolyhedral group}. In  
\S\ref{ultimo}, we shall make use of the following particular groups:
\begin{itemize}
\item $G_6=\sigma(\widetilde{\cT}\times\widetilde{\cT})$ of order $288$;

\item $G_8=\sigma(\widetilde{\cO}\times\widetilde{\cO})$ of order $1152$;

\item $G_{12}=\sigma(\widetilde{\cI}\times\widetilde{\cI})$ of order $7200$.
\end{itemize}

\section{Graph-theoretic intuition}\label{graph}

\subsection{McKay quivers}

If $\fG\subset \SL(n,\IC)$ is a finite subgroup, it defines a natural faithful
representation $\cQ$ of $\fG$. Let $\{V_0,\ldots,V_k\}$ be a complete set of
irreducible representations of $\fG$, where $V_0$ denotes the trivial one. For
each such representation, one may decompose the tensor products
$$
\cQ\otimes V_i\cong\bigoplus_{j=0}^k V_j^{\oplus a_{i,j}}
$$
for some positive integers $a_{i,j}$. If the character of the representation
$\cQ$  is real-valued, then $a_{i,j}=a_{j,i}$ for all $i,j$. One defines the
\emph{McKay quiver} as the unoriented quiver with vertices
$V_0,V_1,\ldots,V_k$ and $a_{i,j}$ edges between the vertices $V_i$ and $V_j$.
In particular, this quiver may contain some loops. For our purpose, we only
consider the \emph{reduced} McKay quiver with vertices $V_1,\ldots,V_k$ and
one edge between $V_i$ and $V_j$ if $i\neq j$ and $a_{i,j}\neq 0$: this means
that we remove from the McKay quiver the vertex $V_0$, all edges starting
from it, all loops and all multiple edges. When there is an edge joining $V_i$
and $V_j$, the vertices are called \emph{adjacent}.

One may check that all finite subgroups of
$\SL(2)$ or $\SO(3,\IR)$ enter in this context since their natural
representation $\cQ$ is real-valued.

\subsection{McKay quivers for the polyhedral groups}

For each binary polyhedral group $\widetilde{G}\subset \SU(2)$ and its
corresponding polyhedral group $G\subset \SO(3,\IR)$, we draw the reduced McKay
quiver with our conventions. For the binary polyhedral groups, we denote by a
white vertex the pure representations and by a black vertex the binary ones. We
get (see for example \cite{GSV,GNS1,GNS2}) the graphs of figure \ref{fig:intersectgraph}.

\newpage

\begin{center}
\begin{figure}
\caption{Reduced McKay quivers}
\label{fig:intersectgraph}
\begin{rotate}{-90}
$
\begin{array}{|l|r|}
\hline
\hspace{5cm}\text{Dimension }2&\text{Dimension }3\hspace{3.2cm}\\
\hline
\xymatrix@C=1.3cm{
\widetilde{C}_n
&*+0{\bullet}\ar@{-}[r]\save[]+<0cm,-0.3cm>*\txt{$\widetilde{\chi}_1$}\restore
&*+0{\circ}\save[]+<0cm,-0.3cm>*\txt{$\chi_1$}\restore\ar@{.}[r]
&*+0{\bullet}\save[]+<0cm,-0.3cm>*\txt{$\widetilde{\chi}_{n-1}$}\restore\ar@{-}[r]
&*+0{\circ}\save[]+<0cm,-0.3cm>*\txt{$\chi_{n-1}$}\restore\ar@{-}[r]
&*+0{\bullet}\save[]+<0cm,-0.3cm>*\txt{$\widetilde{\chi}_n$}\restore
}
&
\xymatrix@C=1.3cm{
*+0{\circ}\ar@{-}[r]\save[]+<0cm,-0.3cm>*\txt{$\chi_1$}\restore
&*+0{\circ}\save[]+<0cm,-0.3cm>*\txt{$\chi_2$}\restore\ar@{.}[r]
&*+0{\circ}\save[]+<0cm,-0.3cm>*\txt{$\chi_{n-2}$}\restore\ar@{-}[r]
&*+0{\circ}\save[]+<0cm,-0.3cm>*\txt{$\chi_{n-1}$}\restore
&C_n
}\\
\hline
\xymatrix@C=1.3cm@R=0.3cm{
&&&&&&*+0{\bullet}\save[]+<0.3cm,0cm>*\txt{$\widetilde{\chi}_1$}\restore\\
\widetilde{D}_{2\ell+1}
&*+0{\circ}\ar@{-}[r]\save[]+<0cm,-0.3cm>*\txt{$\chi_1$}\restore
&*+0{\bullet}\ar@{-}[r]\save[]+<0cm,-0.3cm>*\txt{$\widetilde{\sigma}_1$}\restore
&*+0{\circ}\ar@{.}[r]\save[]+<0cm,-0.3cm>*\txt{$\tau_1$}\restore
&*+0{\bullet}\ar@{-}[r]\save[]+<0cm,-0.3cm>*\txt{$\widetilde{\sigma}_\ell$}\restore
&*+0{\circ}\ar@{-}[ru]\ar@{-}[rd]\save[]+<0cm,-0.3cm>*\txt{$\tau_\ell$}\restore\\
&&&&&&*+0{\bullet}\save[]+<0.3cm,0cm>*\txt{$\widetilde{\chi}_2$}\restore
}
&
\xymatrix@C=1.3cm@R=0.3cm{
\\
*+0{\circ}\ar@{-}[r]\save[]+<0cm,-0.3cm>*\txt{$\chi_1$}\restore
&*+0{\circ}\save[]+<0cm,-0.3cm>*\txt{$\tau_1$}\restore\ar@{.}[r]
&*+0{\circ}\save[]+<0cm,-0.3cm>*\txt{$\tau_{\ell-1}$}\restore\ar@{-}[r]
&*+0{\circ}\save[]+<0cm,-0.3cm>*\txt{$\tau_\ell$}\restore
&D_{2\ell+1}
}\\
\hline
\xymatrix@C=1.3cm@R=0.3cm{
&&&&&&&*+0{\circ}\save[]+<0.3cm,0cm>*\txt{$\chi_2$}\restore\\
\widetilde{D}_{2\ell}
&*+0{\circ}\ar@{-}[r]\save[]+<0cm,-0.3cm>*\txt{$\chi_1$}\restore
&*+0{\bullet}\ar@{-}[r]\save[]+<0cm,-0.3cm>*\txt{$\widetilde{\sigma}_1$}\restore
&*+0{\circ}\ar@{-}[r]\save[]+<0cm,-0.3cm>*\txt{$\tau_1$}\restore
&*+0{\bullet}\ar@{.}[r]\save[]+<0cm,-0.3cm>*\txt{$\widetilde{\sigma}_2$}\restore
&*+0{\circ}\ar@{-}[r]\save[]+<0cm,-0.3cm>*\txt{$\tau_{\ell-1}$}\restore
&*+0{\bullet}\ar@{-}[ru]\ar@{-}[rd]\save[]+<0cm,-0.3cm>*\txt{$\widetilde{\sigma}_\ell$}\restore\\
&&&&&&&*+0{\circ}\save[]+<0.3cm,0cm>*\txt{$\chi_3$}\restore
}
&
\xymatrix@C=1.3cm@R=0.3cm{
&&&&*+0{\circ}\ar@{-}[dd]\save[]+<0.3cm,0cm>*\txt{$\chi_2$}\restore\\
*+0{\circ}\ar@{-}[r]\save[]+<0cm,-0.3cm>*\txt{$\chi_1$}\restore
&*+0{\circ}\ar@{.}[r]\save[]+<0cm,-0.3cm>*\txt{$\tau_1$}\restore
&*+0{\circ}\ar@{-}[r]\save[]+<0cm,-0.3cm>*\txt{$\tau_{\ell-2}$}\restore
&*+0{\circ}\ar@{-}[ur]\ar@{-}[rd]\save[]+<0cm,-0.3cm>*\txt{$\tau_{\ell-1}$} \restore
&&D_{2\ell}\\
&&&&*+0{\circ}\save[]+<0.3cm,0cm>*\txt{$\chi_3$}\restore
}\\
\hline
\xymatrix@C=1.3cm@R=0.3cm{
&&&*+0{\bullet}\ar@{-}[d]\save[]+<0.3cm,0cm>*\txt{$\widetilde{\chi}_1$}\restore\\
\widetilde{\cT}
&*+0{\circ}\ar@{-}[r]\save[]+<0cm,-0.3cm>*\txt{$\chi_1$}\restore
&*+0{\bullet}\ar@{-}[r]\save[]+<0cm,-0.3cm>*\txt{$\widetilde{\chi}_2$}\restore
&*+0{\circ}\ar@{-}[r]\save[]+<0cm,-0.3cm>*\txt{$\chi_3$}\restore
&*+0{\bullet}\ar@{-}[r]\save[]+<0cm,-0.3cm>*\txt{$\widetilde{\chi}_3$}\restore
&*+0{\circ}\save[]+<0cm,-0.3cm>*\txt{$\chi_2$}\restore
}
&
\xymatrix@C=1.3cm@R=0.3cm{
\\
*+0{\circ}\ar@{-}[r]\save[]+<0cm,-0.3cm>*\txt{$\chi_1$}\restore
&*+0{\circ}\ar@{-}[r]\save[]+<0cm,-0.3cm>*\txt{$\chi_3$}\restore
&*+0{\circ}\save[]+<0cm,-0.3cm>*\txt{$\chi_2$}\restore
&& \cT
}
\\
\hline
\xymatrix@C=1.3cm@R=0.3cm{
&&&*+0{\circ}\ar@{-}[d]\save[]+<0.3cm,0cm>*\txt{$\chi_2$}\restore\\
\widetilde{\cO}
&*+0{\bullet}\ar@{-}[r]\save[]+<0cm,-0.3cm>*\txt{$\widetilde{\chi}_1$}\restore
&*+0{\circ}\ar@{-}[r]\save[]+<0cm,-0.3cm>*\txt{$\chi_3$}\restore
&*+0{\bullet}\ar@{-}[r]\save[]+<0cm,-0.3cm>*\txt{$\widetilde{\chi}_3$}\restore
&*+0{\circ}\ar@{-}[r]\save[]+<0cm,-0.3cm>*\txt{$\chi_4$}\restore
&*+0{\bullet}\ar@{-}[r]\save[]+<0cm,-0.3cm>*\txt{$\widetilde{\chi}_2$}\restore
&*+0{\circ}\save[]+<0cm,-0.3cm>*\txt{$\chi_1$}\restore
}
&
\xymatrix@C=1.3cm@R=0.3cm{
*+0{\circ}\ar@{-}[rd]\ar@{-}[dd]\save[]+<-0.3cm,0cm>*\txt{$\chi_2$}\restore\\
&*+0{\circ}\ar@{-}[r]\save[]+<0cm,-0.3cm>*\txt{$\chi_4$}\restore
&*+0{\circ}\save[]+<0cm,-0.3cm>*\txt{$\chi_1$}\restore
&&\cO\\
*+0{\circ}\ar@{-}[ur]\save[]+<-0.3cm,0cm>*\txt{$\chi_3$}\restore
}
\\
\hline
\xymatrix@C=1.3cm@R=0.3cm{
&&&*+0{\circ}\ar@{-}[d]\save[]+<0.3cm,0cm>*\txt{$\chi_2$}\restore\\
\widetilde{\cI}
&*+0{\bullet}\ar@{-}[r]\save[]+<0cm,-0.3cm>*\txt{$\widetilde{\chi}_2$}\restore
&*+0{\circ}\ar@{-}[r]\save[]+<0cm,-0.3cm>*\txt{$\chi_3$}\restore
&*+0{\bullet}\ar@{-}[r]\save[]+<0cm,-0.3cm>*\txt{$\widetilde{\chi}_4$}\restore
&*+0{\circ}\ar@{-}[r]\save[]+<0cm,-0.3cm>*\txt{$\chi_4$}\restore
&*+0{\bullet}\ar@{-}[r]\save[]+<0cm,-0.3cm>*\txt{$\widetilde{\chi}_3$}\restore
&*+0{\circ}\ar@{-}[r]\save[]+<0cm,-0.3cm>*\txt{$\chi_1$}\restore
&*+0{\bullet}\save[]+<0cm,-0.3cm>*\txt{$\widetilde{\chi}_1$}\restore
}
&
\xymatrix@C=1.3cm@R=0.3cm{
*+0{\circ}\ar@{-}[rd]\ar@{-}[dd]\save[]+<-0.3cm,0cm>*\txt{$\chi_2$}\restore\\
&*+0{\circ}\ar@{-}[r]\save[]+<0cm,-0.3cm>*\txt{$\chi_4$}\restore
&*+0{\circ}\save[]+<0cm,-0.3cm>*\txt{$\chi_1$}\restore
&&\cI\\
*+0{\circ}\ar@{-}[ur]\save[]+<-0.3cm,0cm>*\txt{$\chi_3$}\restore
}
\\
\hline
\end{array}
$
\end{rotate}
\end{figure}
\end{center}

\newpage

In the sequel, we shall interpret these graphs as the intersection graphs of a
family of smooth rational curves meeting transversally. One may
then get the following intuition: looking at the two-dimensional graphs, if one
contracts the curves associated to a binary representation (black nodes), then
one gets as intersection graph precisely the corresponding graph in dimension
three!

Another property of the two-dimensional quivers is that no two pure
representations  and no two binary representations are adjacent. This means
that the preceding idea of contraction contracts only one curve each time.

\section{Exceptional fibres in dimensions two and three}\label{exceptional}

Considering the Hilbert-Chow morphism $\cH:\Hilb{\fG}{}{\IC^n}\longrightarrow
\IC^n/\fG$,  our purpose is to describe the exceptional fibre $\cH^{-1}(O)$
over the origin $O~\in~\IC^n/\fG$ in the two- and three-dimensional cases. Note
that all finite subgroups of $\SL(2,\IC)$ or $\SO(3,\IR)$ enter in the context
of \S\ref{s:clusters} since they are subgroups of index $2$ of a
reflection group (see \cite[\S 2.7]{GNS2}). Hence we may apply the general
procedure for the study of the clusters supported at the origin.

The understanding of the exceptional fibre in these cases was achieved by
Ito-Nakamura \cite{INm1,INm2} in dimension two and by
Gomi-Nakamura-Shinoda \cite{GNS1,GNS2} in dimension three, by a case-by-case
analysis. For the two-dimensional case, there is another proof by
Crawley-Boevey \cite{CB} avoiding this case-by-case analysis. We recall 
the results.

For any finite group $\fG$, $\Irr^*(\fG)$ denotes the set of irreducible
representations but the trivial one.  

\subsection{Structure of the exceptional fibre in dimension two}

Let $\widetilde{G}\subset \SL(2,\IC)$ be a binary polyhedral group and denote
the Hilbert-Chow morphism by
$$
\widetilde{\pi}:\Hilb{\widetilde{G}}{}{\IC^2}\longrightarrow \IC^2/\widetilde{G}.
$$
For each non-trivial irreducible representation $\rho$ of $\widetilde{G}$, set
$$
E(\rho):=\{I\in \widetilde{\pi}^{-1}(O)_{\text{red}}\,|\, V(I)\supset \rho\}.
$$

\begin{theorem}\emph{(\cite[Theorem 3.1]{INm1}}
\begin{itemize}
\item Each $E(\rho)$ is a smooth rational curve of self-intersection $-2$.

\item $\widetilde{\pi}^{-1}(O)_{\text{red}}=\bigcup_{\rho} E(\rho)$ and
$\widetilde{\pi}^{-1}(O)=\sum_{\rho} \dim \rho\cdot E(\rho)$ as a Cartier-divisor, $\rho\in\Irr^*(\widetilde{G})$.

\item If $I\in E(\rho)$ and $I\notin E(\rho')$ for all $\rho\neq \rho'$, then
$V(I)\cong \rho$.

\item If $I\subset E(\rho)\cap E(\rho')$, then $V(I)\cong \rho\oplus\rho'$ and
the curves $E(\rho)$ and $E(\rho')$ intersect transversally at $I$.

\item The intersection graph of these curves is the reduced McKay quiver of
the group $\widetilde{G}$.
\end{itemize}
\end{theorem}

In particular, a generator $V(I)$ does not contain more than one copy of any
irreducible representation, and $E(\rho)\cap E(\rho')\neq \emptyset$ if and
only if the representations $\rho$ and $\rho'$ are adjacent.

\subsection{Structure of the exceptional fibre in dimension three}

Let $G\subset \SO(3,\IR)$ be a polyhedral group and denote the Hilbert-Chow
morphism by
$$
\pi:\Hilb{G}{}{\IC^3}\longrightarrow \IC^3/G.
$$
For each non-trivial irreducible representation $\rho$ of $G$, set
$$
C(\rho):=\{J\in \pi^{-1}(O)_{\text{red}}\,|\, V(J)\supset \rho\}.
$$

\begin{theorem}\emph{(\cite[Theorem 3.1]{GNS2})}
\begin{itemize}
\item Each $C(\rho)$ is a smooth rational curve.

\item $\pi^{-1}(O)_{\text{red}}=\bigcup_{\rho} C(\rho)$, $\rho\in\Irr^*(G)$.

\item If $J\in C(\rho)$ and $J\notin C(\rho')$ for all $\rho\neq \rho'$,  then
$V(J)\cong \rho$.

\item The intersection graph of these curves is the reduced McKay quiver of
the group $G$.
\end{itemize}
\end{theorem}

\subsection{Explicit parameterizations}
\label{ss:explicitparam}

Let us explain briefly the explicit parameterizations of the exceptional
curves obtained in \loccit This description holds both
in dimensions two and three so we do it with our general notations. The example of the cyclic group is treated in \S\ref{s:examplecyclic}. As we explained in \S\ref{s:clusters},
$$
\fm_S/\fn_\fG\cong \bigoplus_{\substack{\rho\in \Irr(\fG)\\ \rho\neq
\rho_0}}\bigoplus_{i=1}^{2\dim\rho}V^{(i)}(\rho)
$$
where $\rho_0$ denotes the trivial representation. Thanks to
the exact sequence
$$
0\longrightarrow \fI/\fn_\fG \longrightarrow \fm_S/\fn_\fG \longrightarrow
\fm_S/\fI \longrightarrow 0,
$$
if one wants to parameterize a flat family of clusters over $\IP_1$, one has to
choose, in the trivial sheaf:
$$
\cO_{\IP^1}\otimes \bigoplus_{\substack{\rho\in \Irr(\fG)\\ \rho\neq
\rho_0}}\bigoplus_{i=1}^{2\dim\rho}V^{(i)}(\rho),
$$
a locally free $\fG$-equivariant sheaf affording the regular
representation on each fibre whose quotient is also locally free. The parameterizations are then produced as
follows: one chooses \emph{one} non trivial subbundle
$$
\cO_{\IP_1}(-1)\otimes \rho\hookrightarrow
\cO_{\IP_1}\otimes(V^{(i)}(\rho)\oplus V^{(j)}(\rho))
$$
for some appropriate choice of the indices, and shows that this gives the
required family whose points $\fI$ are characterized by their generator
$$
V(\fI)\subset \IP(V^{(i)}(\rho)\oplus V^{(j)}(\rho)).
$$
That is: once one choice has been made, the other choices are
automatic, and we shall see that they always correspond to a trivial subbundle
(see \ref{rk:decomp}).

\section{Geometric construction}
\label{s:geomconstr}

Let $\widetilde{G}$ be a binary polyhedral group acting on $A=\IC[x,y]$.  Set
$\tau:=\left<\pm 1\right>\subset \widetilde{G}$ and $G:=\widetilde{G}/\tau$ the
associated polyhedral group as before. It is important for the sequel to begin
so, and not to choose the group $G$ with its action on some coordinates first,
as we shall see. We aim to define a regular map
$$
\cS:\Hilb{\widetilde{G}}{}{\IC^2}\longrightarrow \Hilb{G}{}{\IC^3}
$$
inducing a map between the exceptional fibres over the origin.

Since $A^\tau=\IC[x^2,y^2,xy]$, we consider the following composition of ring
morphisms, with $B=\IC[a,b,c]$:
\begin{equation}\label{def:sigma}
\sigma:\xymatrix@1{B\ar@{->>}[r]& B\left/\langle ab-c^2\rangle\right.\ar[r]^-\sim & A^\tau\ar@{->}[r]& A}
\end{equation}
where the identification is defined by $a=x^2,b=y^2,c=xy$. The action  of
$\widetilde{G}$ on $A$ induces an action of $G$ on $A^\tau$. Using the
identification, we can define an action of $G$ on the coordinates $a,b,c$,
inducing an action on $B$ with the property that the cone $K=\langle
ab-c^2\rangle$ is $G$-invariant. This is the reason why we did not fix the
action of $G$ at first: an other choice of identification would induce
an other action of $G$.

Let $I$ be an ideal of $A$ and $J:=\sigma^{-1}(I)$ the corresponding ideal of
$B$.  Observe the following property of the map $\sigma$:
\begin{lemma}\label{bendef}
If $I$ is a $\widetilde{G}$-cluster in $A$, then $J$ is a $G$-cluster in $B$.
Furthermore, if $I$ is supported at the origin, then so is $J$.
\end{lemma}

\begin{proof}
If $I$ is a $\widetilde{G}$-cluster, then $A/I\cong \IC[\widetilde{G}]$. Since
the group $\tau$ is finite, we have
isomorphisms:
$$
B/J\cong A^\tau/I^\tau\cong (A/I)^\tau\cong \IC[\widetilde{G}]^\tau\cong \IC[G],
$$
hence $J$ is a $G$-cluster in $B$. Furthermore, note that $\sigma^{-1}\fm_A=\fm_B$
hence if $I$ is a $\widetilde{G}$-cluster supported at the origin, one has
$I\subset \fm_A$ and then $J\subset \fm_B$, which implies that $J$ is also
supported at the origin (see \S\ref{ss:generalsetup}).
\end{proof}

Therefore, this construction defines set-theoretically a map between the two
moduli spaces of clusters $\cS:\Hilb{\widetilde{G}}{}{\IC^2}\longrightarrow
\Hilb{G}{}{\IC^3}$ by $\cS(I)\overset{\text{Def}}{=}J$. It remains to see that
this map is a regular morphism.

\begin{proposition}\label{regpro}
The map $\cS$ is regular, projective, and induces a map between the exceptional fibres.
\end{proposition}

\begin{proof}\text{}

$\diamond$ In order to get that the map $\cS$ is regular, we show
that it is induced by a natural transformation between the two functors of
points
$$
\Hilbf{\widetilde{G}}{}{\IC^2}(\cdot)\Longrightarrow \Hilbf{G}{}{\IC^3}(\cdot).
$$
Let $T$ be a scheme and $Z\in \Hilbf{\widetilde{G}}{}{\IC^2}(T)$. Then
$Z\subset T\times\IC^2$ is a flat family of $\widetilde{G}$-clusters over $T$
and the map $Z\hookrightarrow T\times\IC^2$ is $\tau$-equivariant (for a
trivial action on $T$). It induces a family
$$
Z/\tau\hookrightarrow T\times(\IC^2/\tau)\hookrightarrow T\times\IC^3
$$
where the quotient $\IC^2/\tau$ is considered as the cone $\langle ab-c^2\rangle$ in
$\IC^3$.  If $T$ is a point, this is precisely our set-theoretic construction
since then if $Z$ is given by an ideal $I$, $Z/\tau$ is given by the ideal
$I^\tau$.

In order to show that $Z/\tau\in \Hilbf{G}{}{\IC^3}(T)$, we have to prove that
this family is flat over $T$. Since this problem is local in $T$, we may assume
that $T$ is an affine scheme, say $T=\Spec R$. Then the family $Z$ is given by
a $\tau$-equivariant quotient $R\otimes A\twoheadrightarrow Q$ so that the
composition $R\hookrightarrow R\otimes_\IC A\twoheadrightarrow Q$ makes $Q$ a
flat $R$-module. The family $Z/\tau$ is then given by the quotient
$$
R\hookrightarrow R\otimes_\IC B \twoheadrightarrow R\otimes_\IC A^\tau
\twoheadrightarrow Q^\tau,
$$
where the quotient $R\otimes_\IC B \twoheadrightarrow R\otimes_\IC A^\tau$  is
induced by tensorization of the quotient $B\twoheadrightarrow A^\tau$. We have
to show that this makes $Q^\tau$ a flat $R$-module. By hypothesis, the functor
$Q\otimes_R-$ in the category of $R$-modules is exact. Since $\tau$ is finite,
the functor $(-)^\tau$ is also exact in this category, and
we note that the functor $Q^\tau\otimes_R-$ is the composition of this two
functors since
$$
Q^\tau\otimes_R N=\left(Q\otimes_R N\right)^\tau
$$
for any $R$-module $N$. Hence the functor $Q^\tau\otimes_R-$ is exact, which
means that the family is flat.

$\diamond$ The composition of ring morphisms (\ref{def:sigma}) gives an
equivariant ring morphism
$$
\xymatrix{\IC[a,b,c]\ar[r]^\sigma\ar@(ld,rd)[]|{G}& \IC[x,y]\ar@(ld,rd)[]|{\widetilde{G}}}
$$
inducing a surjective map at the level of the invariants: 
$\xymatrix{\IC[a,b,c]^G\ar@{->>}[r]& \IC[x,y]^{\widetilde{G}}}$,
hence a closed immersion
$$
\eta:\IC^2/\widetilde{G}=\Spec \IC[x,y]^{\widetilde{G}} \longrightarrow \Spec \IC[a,b,c]^G =\IC^3/G.
$$
Taking more care of the cone $K=\IC^2/\tau$ (in the notations of the introduction), the equivariant map
$$
\xymatrix{\IC^2\ar[r]\ar@(ld,rd)[]|{\widetilde{G}}& \IC^2/\tau\ar[r]\ar@(ld,rd)[]|{G}& \IC^3\ar@(ld,rd)[]|{G}}
$$
induces the $\eta$ map between the quotients:
$$
\eta:\xymatrix@1{\IC^2/\widetilde{G}\ar[r]^-{\sim}& \left. \left(\IC^2/\tau\right)\right/G\ar[r] & \IC^3/G}
$$
sending the origin $O\in \IC^2/\widetilde{G}$ to the origin $O\in\IC^3/G$  and
by definition of $\cS$ the following diagram is commutative:
$$
\xymatrix{\Hilb{\widetilde{G}}{}{\IC^2}\ar[r]^\cS\ar[d]_{\widetilde{\pi}} & \Hilb{G}{}{\IC^3}\ar[d]^\pi\\\IC^2/\widetilde{G}\ar[r]^\eta & \IC^3/G}
$$
This implies that $\cS$ induces a map between the exceptional fibres
$$
\widetilde{\pi}^{-1}(O)\xrightarrow{\cS} \pi^{-1}(O).
$$

$\diamond$ We prove that the map $\cS$ is proper by applying the valuative
criterion of properness. Let $K$ be any field over $\IC$ and $R\subset K$ any valuation ring with
quotient field $K$. Consider a commutative diagram: 
$$
\xymatrix{\Spec K \ar[r]^-\phi\ar[d]^i & \Hilb{\widetilde{G}}{}{\IC^2}\ar[d]^\cS \\
\Spec R\ar[r]^-\psi & \Hilb{G}{}{\IC^3}}
$$
We have to show that there exists a unique factorization
$$
\xymatrix{\Spec K \ar[r]^-\phi\ar[d]^i & \Hilb{\widetilde{G}}{}{\IC^2}\ar[d]^\cS \\
\Spec R\ar[r]^-\psi \ar[ur]^-{\tilde{\phi}}& \Hilb{G}{}{\IC^3}}
$$
making the whole diagram commute.

By modular interpretation, the data of the map $\phi$ consists in an ideal
$I\subset K[x,y]$ such that $K[x,y]/I\cong \IC[\widetilde{G}]\otimes_{\IC}K$
and $K[x,y]/I$ is $K$-flat (it is here trivial since $K$ is a field). Similarly,
the data of the map $\psi$ consists in an ideal $J\subset R[a,b,c]$ such that
$R[a,b,c]/J\cong \IC[G]\otimes_{\IC}R$ and $R[a,b,c]/J$ is $R$-flat. The
commutativity $\cS\circ \phi=\psi\circ i:\Spec K\rightarrow \Hilb{G}{}{\IC^3}$
means the following. Consider the diagram of ring morphisms induced by natural
extension of scalars and base-change from the map $\sigma$:
$$
\xymatrix{R[a,b,c]\ar[r]^-{\sigma_R}\ar@{^(->}[d]&R[x,y] \ar@{^(->}[d]\\
K[a,b,c]\ar[r]^-{\sigma_K} & K[x,y]}
$$
Then the commutativity condition means that $\sigma_K^{-1}(I)=J\cdot K[a,b,c]$.

We are looking for a map $\tilde{\phi}$ such that $\tilde{\phi}\circ i=\phi$
and $\cS\circ \tilde{\phi}=\psi$, \ie for an ideal $\widetilde{I}\subset
R[x,y]$ such that $R[x,y]/\widetilde{I}\cong \IC[\widetilde{G}]\otimes_{\IC}R$
and $R[x,y]/\widetilde{I}$ is $R$-flat, satisfying the conditions
$\widetilde{I}\cdot K[x,y]=I$ and $\sigma_R^{-1}(\widetilde{I})=J$.

A natural candidate is $\widetilde{I}\overset{\text{Def}}{=}I\cap R[x,y]$.  We
have to prove that it satisfies all the conditions and that it is unique for
these properties. Denote by $\nu:K-\{0\}\rightarrow H$ the valuation with
values in a totally ordered group $H$, satisfying the properties:
$$
\nu(x\cdot y)=\nu(x)+\nu(y) \text{ and } \nu(x+y)\geq \min(\nu(x),\nu(y)) \quad \text{for } x,y\in K-\{0\}
$$
and such that $R=\{x\in K\,|\,\nu(x)\geq 0\}\cup\{0\}$. Recall that $R$ is by
definition integral and that a $R$-module is flat is and only if it is
torsion-free (see for instance \cite{AMD,H}).

\begin{itemize}
\item It is already clear that $\widetilde{I}\cdot K[x,y]\subset I$.
Conversely,  Let $P=\sum_{i,j}p_{i,j}x^iy^j\in I$ and $p\in \{p_{i,j}\}$ an
element of minimal valuation. If $\nu(p)\geq 0$, then $P\in \widetilde{I}$.
Else all coefficients of $p^{-1}P$ have positive valuation and so $p^{-1}P\in
\widetilde{I}$. So $P=p\cdot(p^{-1}P)\in \widetilde{I}\cdot K[x,y]$, hence the
equality.

\item By commutativity of the above diagram, 
\begin{align*}
\sigma^{-1}_R(\widetilde{I})&=\sigma^{-1}_R(I\cap R[x,y])\\
&=\sigma^{-1}_K(I)\cap R[a,b,c]\\
&=(J\cdot K[a,b,c])\cap R[a,b,c].
\end{align*}
Is is already clear that $J\subset (J\cdot K[a,b,c])\cap R[a,b,c]$. Conversely,
let $P\in~ (J\cdot K[a,b,c])\cap R[a,b,c]$, decomposed as $P=\sum_\ell
U_\ell\cdot V_\ell$ with $U_\ell\in J$ and $V_\ell\in K[a,b,c]$. As before,
there exists a coefficient $q$ in all $V_\ell$'s of minimal valuation, and we
assume $\nu(q)<0$ (else there is no problem). Then $q^{-1}P\in J$. By
assumption, the $R$-module $R[a,b,c]/J$ is torsion-free, so the multiplication
by $q^{-1}\in R$ is injective. This means that $P\in J$.

\item By definition, we have an $R$-linear inclusion $R[x,y]/\widetilde{I}\hookrightarrow K[x,y]/I$, which shows that $R[x,y]/\widetilde{I}$ is torsion-free, hence flat. It inherits an action of $\widetilde{G}$ and since $K[x,y]/I\cong \IC[\widetilde{G}]\otimes_{\IC}K$, there exists a subrepresentation $V$ of $\IC[\widetilde{G}]$ such that $R[x,y]/\widetilde{I}\cong V\otimes_\IC R$ (this uses the flatness, see \cite[lemma 9.4]{INm2}). By the isomorphism of $R$-modules $R[x,y]/\widetilde{I}\otimes_R K\cong K[x,y]/I$, the representation $V$ is such that $V\otimes_R K=\IC[\widetilde{G}]\otimes_\IC K$, which forces $V\cong\IC[\widetilde{G}]$.

\item The uniqueness of the candidate follows from the condition $\tilde{I}\cdot K[x,y]=I$ since as we already
noted:
$$
I\cap R[x,y]=(\tilde{I}\cdot K[x,y])\cap R[x,y]=\tilde{I}
$$
so our natural candidate is the only possibility.
\end{itemize}

$\diamond$ To finish with, remark that any proper map between to
quasi-projective varieties is automatically a projective map.
\end{proof}

\section{Contracted versus non-contracted fibres}\label{sezteoremone}

\begin{theorem}\label{th:theoremone}
Consider the restriction of the map
$\cS:\Hilb{\widetilde{G}}{}{\IC^2}\longrightarrow \Hilb{G}{}{\IC^3}$ to a
reduced curve $E(\rho)$. Then:
\begin{enumerate}
\item If the representation $\rho$ is pure, then $\cS$ maps isomorphically  the
curve $E(\rho)$ onto the curve $C(\rho)$.

\item If the representation $\rho$ is binary, then $\cS$ contracts the curve
$E(\rho)$ to a point.
\end{enumerate}
\end{theorem}

\begin{proof}
Let $E(\rho)$ be any exceptional curve. Since the map $\cS$ sends this curve to
the bunch of curves $\pi^{-1}(O)$, the image lies in some irreducible component
$C$ and the restricted morphism $\cS:E(\rho)\rightarrow C$ is a proper map. We
prove that:
\begin{itemize}
\item if the representation $\rho$ is binary, then the map
$\cS:E(\rho)\rightarrow C$ contracts the curve to a point;

\item if the representation $\rho$ is pure, then $C=C(\rho)$ and the restricted
map $\cS:E(\rho)\rightarrow C(\rho)$ is an isomorphism.
\end{itemize}

The parameterizations of the two curves $E(\rho)$ and $C$ defines a composite
proper map $f$ whose properties reflect those of the restriction of $\cS$:
$$
\xymatrix{\IP_1\ar[r]^-\sim_-\phi\ar[d]^f&E(\rho)\subset
\Hilb{\widetilde{G}}{}{\IC^2}
\ar[d]^\cS\\
\IP_1\ar[r]^-\sim_-\psi&C\subset\Hilb{G}{}{\IC^3}}
$$
We know (see \cite[II.6.8,II.6.9]{H}) that either the map $f$ contracts the
curve to a point, or is a finite surjective map. The basic idea in order to
determine which case occurs, is to suppose given an ample line bundle
$\cO_{\IP_1}(a)$ on the target (with $a>0$): if the map $f$ contracts the
curve to a point, then $f^*\cO_{\IP_1}(a)$ is trivial and else
$f^*\cO_{\IP_1}(a)\cong \cO_{\IP_1}(\deg(f)\cdot a)$ is ample.

The natural candidate for an ample line bundle over the curve $C$ is the
determinant $\det(p_*\cO_{Z(C)})$ obtained by restriction of the universal family $Z(C):=\left.\cZ_G\right|_C$.

The parameterization $\IP_1\overset{\phi}{\longrightarrow}
\Hilb{\widetilde{G}}{}{\IC^2}$ of the curve $E(\rho)$ corresponds to a flat
family $Z_{\widetilde{G}}(\rho)\subset \IP_1\times\IC^2$ which is the
restriction to $E(\rho)$ of the universal family $\cZ_{\widetilde{G}}$ over
$\Hilb{\widetilde{G}}{}{\IC^2}$. The direct image
$p_*\cO_{Z_{\widetilde{G}}(\rho)}$ is a vector bundle of rank
$|\widetilde{G}|$ over $\IP_1$ equipped with an action of $\widetilde{G}$
affording the regular representation on each fibre. It admits an isotypical
decomposition over the irreducible representation of $\widetilde{G}$ and we
recall the well-known explicit decomposition:
\begin{lemma}\label{lemm:decompfibre}
$$
p_*\cO_{Z_{\widetilde{G}}(\rho)}\cong
\left(\cO_{\IP_1}(1)\oplus\cO_{\IP_1}^{\oplus \dim \rho-1}\right)\otimes
\rho\oplus\bigoplus_{\substack{\rho'\in \Irr(\widetilde{G})\\\rho'\neq
\rho}}\cO_{\IP_1}^{\oplus\dim \rho'}\otimes \rho'
$$
\end{lemma}

\begin{proof}[Proof of the lemma] This is an equivalent form of \cite[\S 2.1 lemma]{KV} or
\cite[Proposition 6.2(3)]{INj}. We recall briefly the argument. 
Since this bundle is a quotient of $\cO_{\IP_1}\otimes A$ (see
\S\ref{ss:explicitparam}), it is generated by its global sections, hence is a
sum of line bundles $\cO_{\IP_1}(a)$ for $a\geq 0$. By the classical observation
$\deg(p_*\cO_{Z_{\widetilde{G}}(\rho)})=1$ (see \cite{GSV}), all line
bundles are trivial but one, of degree one.
\end{proof}

In particular, note that $\det(p_*\cO_{Z_{\widetilde{G}}(\rho)})\cong
\cO_{\IP_1}(\dim \rho)$ is the ample determinant line bundle in dimension two.

Thanks to the functorial definition of the map $\cS$, the composition
$$
\IP_1\overset{\phi}{\longrightarrow}
\Hilb{\widetilde{G}}{}{\IC^2}\overset{\cS}{\longrightarrow} \Hilb{G}{}{\IC^3}
$$
parameterizes the flat family $Z_{\widetilde{G}}(\rho)/\tau$ whose structural
sheaf is
$\cO_{Z_{\widetilde{G}}(\rho)/\tau}=~\left(\cO_{Z_{\widetilde{G}}(\rho)}\right)^\tau$
and one gets:
$$
f^*(\det(p_*\cO_{Z(C)}))=\det\left((p_*\cO_{Z_{\widetilde{G}}(\rho)})^\tau\right).
$$
Now, as we noticed in \S\ref{ss:representations}, taking the invariants under
$\tau$ keeps invariant the pure representations and kills the binary ones.
Hence:
\begin{itemize}
\item If the representation $\rho$ is binary, then:
$$
\left(p_*\cO_{Z_{\widetilde{G}}(\rho)}\right)^\tau\cong \bigoplus_{\rho'\in
\Irr(G)}\cO_{\IP_1}^{\oplus\dim \rho'}\otimes \rho'
$$
hence $\det(p_*\cO_{Z_{\widetilde{G}}(\rho)})^\tau\cong\cO_{\IP_1}$ is trivial;

\item If the representation $\rho$ is pure, then:
$$
\left(p_*\cO_{Z_{\widetilde{G}}(\rho)}\right)^\tau\cong
\left(\cO_{\IP_1}(1)\oplus\cO_{\IP_1}^{\oplus \dim \rho-1}\right)\otimes
\rho\oplus\bigoplus_{\substack{\rho'\in \Irr(G)\\\rho'\neq
\rho}}\cO_{\IP_1}^{\oplus\dim \rho'}\otimes \rho'
$$
hence $\det(p_*\cO_{Z_{\widetilde{G}}(\rho)})^\tau\cong\cO_{\IP_1}(\dim \rho)$
is ample.
\end{itemize}

This achieves the first part of the proof. It remains to show that in the case
of a pure representation $\rho$, the target curve is $C=C(\rho)$ and that the
finite surjective map $f$ is an isomorphism. We do it by hand. A point $I\in
E(\rho)$ is characterized by the choice of $V(I)$ and generically $V(I)\cong
\rho$. For a pure representation $\rho$, the polynomials defining $V(I)$ are
even hence:
$$
V(I^\tau)=V\left((A\cdot V(I)+\fn_A)^\tau\right)\supset V(I)
$$
so generically $V(I^\tau)=V(I)$ (only modified by setting $a=x^2,b=y^2,c=xy$).
This means that $C=C(\rho)$ and if $I\neq J\in E(\rho)$, then $V(I)\neq V(J)$
hence the images are also different, so the map is generically injective. This
concludes the proof.
\end{proof}

As a byproduct of our argument, we get the following equivalent in dimension three
of the lemma \ref{lemm:decompfibre} which, to our knowledge, does not appear explicitly 
in the literature:

\begin{corollary}
\label{cor:decomp3}
For any finite subgroup $G\subset\SO(3,\IR)$ and any non-trivial representation $\rho$ of $G$, the restriction of the tautological bundle to the exceptional curve $C(\rho)$ decomposes as:
$$
p_*\cO_{Z_{G}(\rho)}\cong
\left(\cO_{\IP_1}(1)\oplus\cO_{\IP_1}^{\oplus \dim \rho-1}\right)\otimes
\rho\oplus\bigoplus_{\substack{\rho'\in \Irr(G)\\\rho'\neq
\rho}}\cO_{\IP_1}^{\oplus\dim \rho'}\otimes \rho'
$$
\end{corollary}

\begin{proof}
The same argument as in the proof of lemma \ref{lemm:decompfibre} shows that this bundle in generated
by its global sections. The bijectivity of the map $f$ on the curves associated to pure
representations (in the notation of the proof of theorem \ref{th:theoremone}) implies that $\det(p_*\cO_{Z_{G}(\rho)})\cong
\cO_{\IP_1}(\dim \rho)$, hence in the isotypical decomposition there is only one non-trivial line bundle, of degree one, and
we already know by the explicit parameterizations that the isotypical component corresponding to $\rho$ is not trivial.
\end{proof}

\begin{remark}
\label{rk:decomp}
In the decomposition of the lemma \ref{lemm:decompfibre}, the unique presence
of the $\cO_{\IP_1}(1)$ corresponds to the choice of the line $V(I)$ in a
projective space $\IP(\rho\oplus \rho)$ as explicitly described in
\S\ref{ss:explicitparam}. The fact that no other ample bundle occurs reflects
the property that once one choice has been made, the other generators of the
ideal do not involve the choice any more, as one can easily notice from the
explicit computations of \cite[\S 13,\S 14]{INm2} (see \S\ref{s:examplecyclic} in this paper
for an example). In the three-dimensional case, the same situation occurs
thanks to the corollary \ref{cor:decomp3}.
\end{remark}

We get now the theorem \ref{th:barth} presented in the introduction as a
corollary of the theorem \ref{th:theoremone}:

\begin{corollary}
The image $\cY:=\cS(\Hilb{\widetilde{G}}{}{\IC^2})$ projects onto the quotient
$K/G$, inducing a partial resolution of singularities containing only the
exceptional curves corresponding to pure representations. The map
$\cS:\Hilb{\widetilde{G}}{}{\IC^2}\longrightarrow \cY$ is a resolution of
singularities contracting the excess exceptional curves to ordinary nodes.
\end{corollary}

\begin{proof}
The projection $\pi:\cY\longrightarrow \IC^3/G$ factors through $K/G$ by
construction of $\cY$. The other assertions result from theorem
\ref{th:theoremone}. The excess curves contract to ordinary nodes since, as one checks with the
figure \ref{fig:intersectgraph}, each excess $(-2)$-curve is contracted to a different point.
\end{proof}

\section{Example: the cyclic group case}
\label{s:examplecyclic}

Let the cyclic group $\widetilde{C}_n\cong \IZ/(2n)\IZ$ act on $\IC^2$ with
generator:
\begin{eqnarray*}\left(\begin{array}{cc}
\xi&0\\
0&\xi ^{-1}\\
\end{array}\right)\quad \text{with }\xi = e^{\frac{2\pi \ci}{(2n)}}.
\end{eqnarray*}
The choice of coordinates made in \S\ref{s:geomconstr} implies that the group
$C_n\cong\IZ/n\IZ$ acts on $\IC^3$ with generator:
$$
\left(\begin{array}{ccc}
\xi ^2&0&0\\
0&\xi ^{-2}&0\\
0&0&1
\end{array}\right).
$$
The irreducible representations of the cyclic group $\widetilde{C}_{n}$ are
given by the matrices $(\xi^i)$, $i=0,\ldots,2n-1$. For $i$ even, they
are also the irreducible representations of $C_n$. There are then $n$ pure
and $n$ binary representations. With the notations of
\S\ref{ss:representations}, we set $\chi_i:=\rho_{2i}$ and
$\widetilde{\chi}_i=\rho_{2i+1}$ for $i=0,\ldots,n-1$. By Theorem
\ref{th:theoremone}, the exceptional curves on $\Hilb{\widetilde{C}_n}{}{\IC^2}$
corresponding to the binary representations are contracted by $\cS$ to a node on
$\cS(\Hilb{\widetilde{C}_n}{}{\IC^2})$ whereas the curves corresponding to the
pure representations are in $1:1$ correspondence with the exceptional curves
downstairs (see figure \ref{fig:contractioncyclic}). In this section, we check this by a direct computation.

\begin{center}
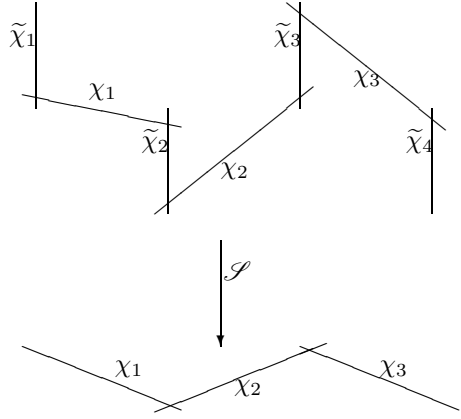
\begin{figure}[h!]
\begin{picture}(170,160)
\put(10,120){\line(0,1){40}} \put(60,80){\line(0,1){40}}
\put(110,120){\line(0,1){40}}\put(160,80){\line(0,1){40}}
\put(5,125){\line(5,-1){60}} \put(55,80){\line(5,4){60}}
\put(105,160){\line(5,-4){60}} \put(0,145){\makebox{$\widetilde{\chi}_1$}}
\put(50,105){\makebox{$\widetilde{\chi}_2$}}\put(100,145){\makebox{$\widetilde{\chi}_3$}}
\put(150,105){\makebox{$\widetilde{\chi}_4$}}\put(30,125){\makebox{$\chi_1$}}
\put(80,95){\makebox{$\chi_2$}} \put(130,130){\makebox{$\chi_3$}}
\put(80,70){\vector(0,-1){40}}\put(81,55){\makebox{$\cS$}}
\put(5,30){\line(5,-2){60}}\put(55,5){\line(5,2){65}}
\put(110,30){\line(5,-2){60}} \put(40,20){\makebox{$\chi_1$}}
\put(85,13){\makebox{$\chi_2$}} \put(140,20){\makebox{$\chi_3$}}
\end{picture}
\caption{Contracted fibres for $\widetilde{C}_4$}
\label{fig:contractioncyclic}
\end{figure}
\end{center}

The ring of invariants
$\IC[x,y]^{\widetilde{C}_n}$ is generated by $x^{2n}$, $y^{2n}$, $xy$ and
$\IC[a,b,c]^{C_n}$ is generated by $c, a^n, b^n, ab$.
Recall the description of the exceptional curves of
$\Hilb{\widetilde{C}_n}{}{\IC^2}$ following \cite[Theorem 12.3]{INm2}. We sort the basis
of the algebra of coinvariants with respect to each
irreducible representation:
$$
\{1\},\{x,y^{2n-1}\},\ldots,\{x^i,y^{2n-i}\},\ldots,\{x^{2n-1},y\}.
$$
To choose a cluster $I/\fn_A$ supported at the origin amounts in choosing
one copy of each non-trivial representation, \ie for all $i=1,\ldots,2n-1$
a point $(p_i:q_i)\in\IP_1$ defining the ideal by the generators:
$$
\langle p_1x-q_1y^{2n-1},\ldots,p_ix^i-q_iy^{2n-i},\ldots,p_{2n-1}x^{2n-1}-q_{2n-1}y\rangle.
$$
But the point is that one only needs \emph{one} choice. Suppose there exists an
index $i$ such that $p_iq_i\neq 0$, and take the smaller $i$ with this
property. Set $p=p_i,q=q_i$ and $v=px^i-qy^{2n-i}$. Then since $xy$ is
invariant, $x^{i+1},\ldots,x^{2n-1}\in I/\fn_A$ and
$y^{2n-i+1},\ldots,y^{2n-1}\in I/\fn_A$ so all our other choices were trivial,
and $V(I)=\IC\cdot v$. More formally, we parameterized the exceptional curve
$E(\rho_i)$ by a subbundle:
$$
\cO_{\IP_1}(-1)\otimes \rho_i\oplus\bigoplus_{j\neq i}\cO_{\IP_1}\otimes
\rho_j\hookrightarrow\bigoplus_{j}(\cO_{\IP_1}\oplus \cO_{\IP_1})\otimes \rho_j.
$$
If there is no such index, suppose $x^i$ is the minimal power of $x$ in the choice: in
order to find once each non-trivial representation one has to choose
$y^{2n-i+1}$ and the minimal set of generators $V(I)=~\IC\cdot x^i\oplus\IC\cdot
y^{2n-i+1}$ contains two adjacent representations.

Otherwise stated, a $\widetilde{C}_n$-cluster at the origin takes the form:
$$
\begin{array}{c}
I_j(p:q):=\langle px^j-qy^{2n-j},xy,x^{j+1},y^{2n-j+1}\rangle,\\
1\leq j\leq 2n-1,~~(p:q)\in\IP_1
\end{array}
$$
(the above expression contains enough generators to include the two possible cases) and
$$
E(\rho_j)=\{I_j(p:q)\}.
$$

By the same method, one sees easily that a $C_n$-cluster at the origin takes
the form:
$$
\begin{array}{c}
J_k(s:t):=\langle sa^k-tb^{n-k},c,a^{k+1},b^{n-k+1},ab\rangle,\\
1\leq k\leq n-1,~~(s:t)\in\IP_1
\end{array}
$$
and
\begin{eqnarray*}
C(\chi_k)=\{J_k(s:t)\}.
\end{eqnarray*}
Recall that with the construction (\ref{def:sigma}) we have to compute
$\sigma^{-1} (I_j(p:q))$. Denoting by $\bar{\sigma}$ the map $ B\left/\langle
ab-c^2\rangle\right.\longrightarrow A$, it is equivalent to compute
$\bar{\sigma}^{-1} (I_j(p:q))$. First we compute $I_j(p:q)^{\tau}\in A^{\tau}$.
We distinguish two cases:
\begin{itemize}
\item $j$ \emph{even}, \ie $j=2j',~j'=1,\ldots,n-1$. In this case we have
\begin{eqnarray*}
I_j(p:q)^{\tau}=I_j(p:q)=\langle p{(x^2)}^{j'}-q{(y^2)}^{n-j'}, xy,
{(x^2)}^{j'+1}, {(y^2)}^{n-j'+1}\rangle
\end{eqnarray*}
expressed in $A^{\tau}=\IC[x^2,y^2,xy]$. Then
\begin{align*}
\bar{\sigma}^{-1}(I_j(p:q))&=\langle pa^{j'}-qb^{n-j'}, c, a^{j'+1}, b^{n-j'+1}\rangle\\
&=J_{j'}(p:q).\\
\end{align*}
\item $j$ \emph{odd}, \ie $j=2j'+1,~j'=0,\ldots,n-1$. Observe that $xy\in
I_j(p:q)^{\tau}$ and $(x^2)^{j'+1}$, $y^{n-j'}\in I_j(p:q)^{\tau}$,  but
$px^{2j'+1}-qy^{2n-2j'-1}\notin I_j(p:q)^{\tau}$. So
$$
\bar{\sigma}^{-1}(I_j(p:q))=\langle a^{j'+1}, b^{n-j'}, c\rangle.\\
$$
We observe then that
$$\bar{\sigma}^{-1}(I_j(p:q))\in C(\rho_{j'})\cap
C(\rho_{j'+1})
$$
since
$$
\bar{\sigma}^{-1}(I_j(p:q))= J_{j'}(0:1)=J_{j'+1}(1:0).
$$
\end{itemize}
The curves  $E(\rho_j)$ with $j$ even correspond to the pure representations
and  are not contracted by $\cS$ as the previous computation shows, the curves
with $j$ odd correspond to the binary representations: these are contracted  by $\cS$.

\section{Application}
\label{ultimo}

\subsection{Pencils of symmetric surfaces}

Let $\IH_\IC:=\IH\otimes_\IR\IC$ be the complexification of the space of quaternions. By the choice of the coordinates
$q=a\cdot\Iun+b\cdot\Ii+c\cdot\Ij+d\cdot\Ik$, $a,b,c,d\in \IC$, one gets an isomorphism $\IP_3\cong \IP(\IH_\IC)$ such, that for $n=6,8,12$ the bipolyhedral group $G_n$ acts linearly on $\IP_3$, leaving invariant the quadratic polynomial $Q:=a^2+b^2+c^2+d^2$.

In \cite{Sa} is shown that the next non-trivial invariant is a homogeneous polynomial $S_n$ of degree $n$. Consider then the following pencil of $G_n$-symmetric surfaces in $\IP_3$:
$$
X_n(\lambda)=\{S_n+\lambda Q^{n/2}=0\}, \quad \lambda\in \IC.
$$
In \cite{Sa} is proved that the general surface $X_n(\lambda)$ is smooth and that for each $n$ there are precisely four singular surfaces in the corresponding pencil: the singularities of these surfaces are ordinary nodes forming one orbit through $G_n$.

Consider now the pencil of quotient surfaces in $\left.\IP_3\right/G_n$:
$$
\{\left.X_n(\lambda)\right/G_n\}, \quad \lambda\in \IC.
$$
In \cite{BaSa} is proved that these quotient surfaces have only A-D-E singularities and that the minimal resolutions of singularities $Y_n(\lambda)\rightarrow \left.X_n(\lambda)\right/G_n$ are K3-surfaces with Picard number greater than $19$. For the four nodal surfaces in each pencil, a carefully study of the stabilizors of the nodes shows that, if $X$ denotes one of these nodal surfaces, the image of the node on $X/G_n\subset \left.\IP_3\right/G_n$ is a particular quotient singularity locally isomorphic to $\IC^2/\widetilde{G}\subset \IC^3/G$ for some polyhedral group $G$ explicitly computed (see \cite[\S 3, Proposition 3.1]{BaSa}):
\begin{itemize}
\item for $n=6$: $C_3,\cT$;

\item for $n=8$: $D_2,D_3,D_4,\cO$;

\item for $n=12$: $D_3,D_5,\cT,\cI$.
\end{itemize}

Therefore, our theorem \ref{th:barth} gives locally a group-theoretic interpretation of the exceptional curves of the K3-surfaces $Y_n(\lambda)$ over the particular singularities of the nodal surfaces.

\nocite{*}
\bibliographystyle{amsplain}
\bibliography{BSbiblio}

\end{document}